\documentclass[reqno]{amsart}
\usepackage{amsmath}
\allowdisplaybreaks[0]
\usepackage{paralist}
\usepackage{graphics}
\usepackage{epsfig}
\usepackage[pdftex]{hyperref}
\usepackage{amsfonts}
\usepackage{CJK}
\usepackage{fancyhdr}
\usepackage{graphicx}
\usepackage[dvipsnames,usenames]{color}
\usepackage{titletoc}
\usepackage{latexsym}
\usepackage{amssymb}
\usepackage{multicol}
\usepackage{graphics}
\usepackage{subfigure}
\usepackage{indentfirst}
\usepackage{cases}
\usepackage{curves}
\usepackage{cite}
\newtheorem{theorem}{Theorem}[section]
\newtheorem{corollary}{Corollary}[section]

\newtheorem{definition}{Definition}[section]

\newtheorem{example}{Example}[section]

\numberwithin{equation}{section}

\def\x#1{(\ref{#1})}
\def\R{{\Bbb R}}

\def\N{{\Bbb N}}

\def\bc{\begin{center}}
\def\ec{\end{center}}
\def\ba{\begin{array}}
\def\ea{\end{array}}
\def\be{\begin{equation}}
\def\ee{\end{equation}}
\def\bea{\begin{eqnarray}}
\def\eea{\end{eqnarray}}
\def\beaa{\begin{eqnarray*}}
\def\eeaa{\end{eqnarray*}}

\def\ben{\begin{enumerate}}
\def\een{\end{enumerate}}
\def\hh{\!\!\!\!}
\def\EM{\hh &   &\hh}
\def\EQ{\hh & = & \hh}

\def\LE{\hh & \le & \hh}
\def\GE{\hh & \ge & \hh}

\def\nn{\nonumber}
\def\oo{\infty}
\def\ifl{\iffalse}
\def\lb{\label}
\def\prf{\mbox{\bf Proof.~}}

\title[]{Stability of Lyapunov Exponents, Weak Integral Separation and Nonuniform   Dichotomy  Spectrum}

\author[H. Zhu, Z. Li, X. He]
{Hailong Zhu$^{1}$,\quad Zhaoxiang Li$^{2}$,  Xiuli He$^{3}$}

\address{$^1$ School of Statistics and Applied Mathematics, Anhui
University of Finance and Economics, Bengbu 233030, China}
\address{$^2$Department of Mathematics, Shanghai Normal University, Shanghai, 200234, China}
\address{$^3$  Department of Mahthematics,  Hohai University, Nanjing, 210098, China}

\email{hai-long-zhu@163.com (H. Zhu)}
\email{zxli@shnu.edu.cn (Z. Li)}
\email{hexiu00@163.com (X. He)}

\subjclass[2000]{34D08, 34D09, 34D20} \keywords{Lyapunov
exponents;
Nonuniform exponential dichotomy spectrum; Weak integral separation.}

\begin{document}

\begin{abstract} In this paper, a necessary and sufficient condition for the stability of Lyapunov exponents of linear differential system  are proved in the sense that  the equations satisfy the weaker form of integral separation instead of its classical one.
Furthermore, the existence of full nonuniform exponential dichotomy
spectrum under the condition of weak integral separateness is also presented.
\end{abstract}

\maketitle

\section{\bf{Introduction}}
\setcounter{equation}{0} \noindent

Lyapunov exponents, or characteristic exponents was originated in the thesis of Lyapunov more than 100 years ago \cite{lv-92}. Since then it has become a very important part of the general theory of dynamical systems, and has played a pivotal role in the study of growth rates of time dependent differential equations. For an $n$-dimensional problem, there are $n$ Lyapunov exponents, and it is natural to think about the stability of the  Lyapunov exponents of an $n$-dimensional system. However,  it is unable to identify the stability of the  Lyapunov exponents of a general system even if  for  a regular system with different Lyapunov exponents.
An example from \cite[p. 171]{ad-95} shows that a two-dimensional system \[
\dot{x}_{1}=(1+\frac{\pi}{2}sin(\pi \sqrt{t}))x_{1},\quad \dot{x}_{2}=0
\]
has two distinct Lyapunov exponents $\lambda_{1}=1$ and $\lambda_{1}=0$. This system is regular but not stable.

This stability theory of Lyapunov exponents has attracted the attention of many leading mathematicians since the birth of the Lyapunov exponents. One of the first sufficient conditions for the stability of Lyapunov exponents of system \be\lb{eq1}\dot{x}=A(t)x\ee under
small perturbations of the coefficient matrix can   be  traced back to Perron \cite{per2}.
After that, important results on stability of Lyapunov exponents were repeatedly improved by Bylov, Vinograd, Izobov, Grobman, Million\v{s}\v{c}ikov and several others \cite{by-54, by-65, by-66,M-69, vi-60}, until the necessary condition is established independently by Bylov et al. \cite{by-69} and by Million\v{s}\v{c}ikov \cite{M-69-2} under a key assumption: integral separation. At this point, a necessary and sufficient condition for the stability of Lyapunov exponents has been established. see e.g, \cite[Chap. 5]{ad-95} and \cite{by-65, by-69}  for details.

Based on this observation, the notion of integral separateness  plays a key role in the theory of dynamical systems. This concept has several important features, of which we mention two. First, an integrally separated  system can be reduced   to a diagonal one \cite{by-65} (see also \cite[Theorem 5.3.1]{ad-95}). Second, the stability of Lyapunov exponents \x{eq1} can be established  if and only if this system are integrally separated and \x{eq1} has different Lyapunov exponents $\lambda_{1} > \cdots > \lambda_{n}$ (see \cite{by-69},  \cite[Thm. 5.4.8]{ad-95} for details).

Here we emphasize that there exist large classes of linear differential equations possessing integral separateness.
Furthermore, the
corresponding theory and its applications are widely developed.
In this respect, we can mention,
for example, the classical series of papers \cite{dv-02, pa-79, pa-82, pa-04}, which in particular discuss that \x{eq1} has a full Sacker-Sell spectrum implies that \x{eq1} are integrally separated \cite{bs-00}. For a detailed discussion and historical comments of this concept, we strongly recommend the book \cite{ad-95}.

In this paper, we propose and discuss the weaker form of integral separation instead of its classical one, since there exists linear system
without the existence of integral separateness even if Lyapunov exponents are a set of $n$ different points. For example, consider the following diagonal system
\be\lb{b10} \dot{X}=\left(
  \begin{array}{cc}
    \omega_{1} & 0\\
    0 & \omega_{2}t\sin t \\
  \end{array}
\right)X,\ee
with $\omega_{1}>\omega_{2}>0$ be real paraments. Example \ref{example1} below shows that \x{b10} is not integrally separated.

The first purpose of this study is to replace the condition of integral separateness with weak integral separateness (see Definition \ref{defn22}), and establish a necessary and sufficient condition for the stability of Lyapunov exponents of \x{eq1}.  Here the stability means that the perturbations of the coefficient matrix  are under exponentially decaying. To do this, we need to
extend the known results of bounded growth (see e.g., \cite{pa-79}) to nonuniformly bounded growth, i.e.,
\x{eq1} has a nonuniformly bounded growth if there exist constants $K>0$, $\tilde{a}>0$ and $\tilde{b}\ge 0$ such that
\begin{equation}\lb{eq2}
\|\Phi(t)\Phi^{-1}(s)\|\le K e^{\tilde{a}|t-s|}e^{\tilde{b} s}, \quad {\rm for} \quad t,s \ge 0,
\end{equation}
where $\Phi(t)$ is a fundamental matrix of \x{eq1}.

Our first main result contained in the following theorem extends the result of \cite{by-69} (see also \cite[Thm 5.4.7]{ad-95}).

\begin{theorem}\lb{main1}Assume that {\rm\x{eq1}} with nonuniformly bounded growth  has distinct Lyapunov exponents $\lambda_{1} > \cdots > \lambda_{n}$.
Then they are stable with perturbations of the coefficient matrix being exponentially decaying, i.e., for a purturbed system
 \be\lb{eq3}\dot{x}=(A(t)+B(t))x\ee
 with $\|B(t)\|\le \delta e ^{-\beta t}$ for some $ \delta >0,~\beta\ge 0$, the Lyapunov exponents of  {\rm\x{eq1}} are stable if and only if there exists a Lyapunov transformation $y\leftarrow T(t)x$ transforming {\rm\x{eq1}} to the diagonal form \be\lb{b6}\dot{y}=diag [a_{1}(t),\ldots, a_{n}(t)]y,\ee with $a_{i}$ are weakly integrally separated functions.
\end{theorem}

In addition, as a corollary of Theorem \ref{main1}, another equivalent condition for the stability of Lyapunov exponents, focusing on the properties of fundamental matrix solution but not on the coefficients of \x{eq1}, is given in Section 2.

The second purpose of this paper is to  establish the connection between weak integral separateness and nonuniform exponential dichotomy spectrum \cite{chu-15, zh-14}. Under the condition of full nonuniform exponential dichotomy spectrum, the existence of weak integral separateness  is obvious, since a full Sacker-Sell spectrum implies that \x{eq1} is integrally separated \cite{bs-00}.   However, the contrary is not true in general.
Example \ref{example31} below illustrates that a linear differential system are weakly integrally separated, which does not have a full nonuniform exponential dichotomy spectrum.
The following theorem establishes the existence of full nonuniform exponential dichotomy spectrum under the condition of weak integral separateness with additional support.

\begin{theorem}\lb{main11}
 Assume that {\rm\x{eq1}} with nonuniformly bounded growth is weakly integrally separated. Considering a fundamental matrix solution $\Phi(t)=(\Phi_{1}(t),\ldots,\Phi_{n}(t))$ of {\rm\x{eq1}} such that the distinct Lyapunov exponents of the columns of $\Phi$ are ordered as $\lambda_{1} > \cdots > \lambda_{n}$. If there exists a Lyapunov transformation $y\leftarrow T(t)x$ transforming {\rm\x{eq1}} into the diagonal form {\rm\x{b6}}. For any interval $(\lambda_{i}, \lambda_{i+1})$, there exists a $\lambda\in (\lambda_{i}, \lambda_{i+1})$ such that
      \[\dot{y}=\left(
       \begin{array}{ccccccc}
         a_{1}(t) &   &  &   &   &   &   \\
           & \ddots &   &  &   &   &  \\
           &   & a_{i}(t)&  &  &   &  \\
          &  &  & \lambda &  &  &  \\
           &  &   &  & a_{i+1}(t) &   &   \\
           &  &  &   &   & \ddots &   \\
           &   &   &   &   &   &  a_{n}(t) \\
       \end{array}
     \right)y
  \]
     is weakly integrally separated. Then the nonuniform dichotomy spectrum $\Sigma_{NED}(A)$ is a full spectrum. This full spectrum is constituted by the  disjoint  union  of $n$ closed intervals, i.e.,
    \[\Sigma_{NED}(A)=[a_{1},b_{1}]\cup [a_{2},b_{2}]\cup \cdots [a_{n-1},b_{n-1}]\cup [a_{n},b_{n}],\]
    where $a_{1}\le b_{1}< a_{2}\le b_{2}<\cdots < a_{n}\le b_{n}$.
\end{theorem}

Furthermore, under the condition of weak integral separateness, a  sufficient condition for the existence of a nonuniform exponential dichotomy  is also given in Section 3.


%
%

\section{\bf{Weak integral separateness and stability of Lyapunov exponents}}
\setcounter{equation}{0} \noindent

One of the instability of Lyapunov exponents of \x{eq1} with small perturbations of the coefficients was first shown by Perron. In \cite{per} Perron uses an example of a two-dimensional system to show that small perturbations of coefficients of a linear system can give rise to large shifts of the Lyapunov exponents. For example, consider  a
2-dimensional system (see \cite[p. 135]{ad-95} for details)
\[ \dot{x}_{1} = 0,\quad \quad \quad
\dot{x}_{2} = \pi \sin \pi \sqrt{t} x_{2},
\]
the perturbation
\[Q(t)=\left(
         \begin{array}{cc}
           0 & \delta/\sqrt{t} \\
           \delta/\sqrt{t} & 0 \\
         \end{array}
       \right)
\]
shifts the greatest Lyapunov exponent of this system
by one to the right even though $\delta$ is  sufficiently small.    To illustrate the stability of the Lyapunov exponents under small perturbation, a kind of  stability of the Lyapunov exponents is established as follows.


\begin{definition}{\rm(see \cite[Def. 5.2.1]{ad-95})} The Lyapunov exponents $\lambda_{1}\ge \cdots\ge \lambda_{n}$ of {\rm \x{eq1}} are said to be stable if for any $\varepsilon>0$ there exists a $\delta>0$ such that $\sup_{t\in \R^{+}}\|B(t)\|<\delta$ implies
\be\lb{b1}|\lambda_{i}-\lambda'_{i}|<\varepsilon,\quad i=1,\ldots, n.\ee
\end{definition}


In the theory of differential equations, integral separation, together with some of its
extensions, and modifications (for example, exponential separation, see, e.g., \cite{by-66, pa-82, pa-04}), plays a major role in the theory of diagonalization and  stability of Lyapunov exponents. On the other hand, the notion of integral separateness demands considerably from the dynamics and it is therefore of increasing interest to look for more general types about integral separation.
Based on this observation, we first introduce the definition of the weak integral separateness, and then present an example which is weekly integrally separated but not integrally separated.

\begin{definition}\lb{defn22}The continuous functions $g_{i}, i=1, \ldots, n$, are said to be weakly
integrally separated if for $i=1, \ldots, n-1$,
there exist some costants $a, b \geq 0$ and $D\in \R$ such that
\[
\int^{t}_{s}(g_{i+1}(\tau)-g_{i}(\tau)) d \tau\geq a(t-s)-b s + D, \quad t\ge s \ge 0.
\]
\end{definition}

This definition mimics the classical notion of integral separateness (see e.g., \cite[Def. 5.3.1]{ad-95} for details). Obviously, integral separateness implies week integral separateness due to the fact $b\ge 0$, but not vice versa. We now present an example of weekly integrally separated that is not integrally separated.

\begin{example}\lb{example1} Let $\omega_{1}>\omega_{2}>0$ be real paraments, the functions $\omega_{1}, \omega_{2}t\sin t$ are not integrally separated but are weekly integrally separated.
\end{example}
\noindent\prf{ It is easy to verify that
\bea\lb{b3}\EM \int_{s}^{t}(\omega_{1} -\omega_{2}\tau\sin \tau)d \tau\nn\\
\EQ \omega_{2}t\cos t  -\omega_{2} s \cos s- \omega_{2}\sin t+ \omega_{2}\sin s+\omega_{1}(t-s)\nn\\
\EQ (\omega_{1}-\omega_{2})(t-s)+ \omega_{2}t(\cos t+1) -\omega_{2}s(\cos s+1) + \omega_{2}(\sin s-\sin t)\nn\\
\GE (\omega_{1}-\omega_{2})(t-s) -2\omega_{2}s - 2\omega_{2}. \eea
Furthermore, if $t=2k\pi+\pi$ and $s=2k\pi$ with $k\in \N$, then
\be\lb{b4}\int_{s}^{t}(\omega_{1} -\omega_{2}\tau\sin \tau)d \tau=(\omega_{1}-\omega_{2})(t-s) -2\omega_{2}s.\ee
Thus, the functions $\omega_{1}, \omega_{2}t\sin t$  admit a week integral separateness. By \x{b4}, the perturbation $-2\omega_{2}s$ in \x{b3} can not be eliminated. This means that the integral separateness is not satisfied.\hspace{\stretch{1}}$\Box$}

The following definition about week integral separateness is introduced for  the system \x{eq1}, which pays more attention to the properties of fundamental matrix solution than the coefficients.

\begin{definition}Considering a fundamental matrix solution $\Phi(t)=(\Phi_{1}(t),\ldots,\Phi_{n}(t))$ of {\rm \x{eq1}}. {\rm \x{eq1}} is said to be weakly integrally separated if for $i=1, \ldots, n-1$, there exist some constants $D>0$ and $a, b \geq 0$ such that
\be\lb{b5}
\frac{\|\Phi_{i+1}(t)\|}{\|\Phi_{i+1}(s)\|} \cdot \frac{\|\Phi_{i}(s)\|}{\|\Phi_{i}(t)\|}\geq D e^{a(t-s)-b s}, \quad t\ge s \ge 0.
\ee
\end{definition}

The following theorem compares two definitions of weak integral separateness for a diagonal system, i.e., the connection between the coefficients and the fundamental matrix .
The conclusion can be  easily proved due to the fact that fundamental matrix solution $\Phi(t)=(\Phi_{1}(t),\ldots,\Phi_{n}(t))$ can be expressed by the diagonal elements,  and thus the proof is omitted. In Corollary \ref{main2}, we will further investigate the relationship between the coefficients and the fundamental matrix solution for the system \x{eq1} with a general form under the condition of Lyapunov transformation.

\begin{theorem}The diagonal system {\rm \x{b6}} is weakly integrally separated if and only if its diagonal coefficients are weakly integrally separated.
\end{theorem}

In \cite{pa-79}, Palmer consider the coefficient matrix $A$ of \x{eq1} in a Banach space $\mathfrak{B}$ with the norm $\|A\|=\sup_{t\ge 0}\|A(t)\|$, and  indicate that integral separateness forms an open and dense subset of $\mathfrak{B}$ (see also \cite{M-69}). Therefore weak integral separateness is a generic property in $\mathfrak{B}$ due to the fact that integral separateness is contained in  weak integral separateness.

Before proceeding further, we recall  some notations and notions, which are the keys to illustrating our main results.

\begin{definition}
  A smooth invertible change of variables $y\leftarrow T^{-1}(t)x$ is called a Lyapunov transformation if $T$, $T^{-1}$, and $\dot{T}$ are bounded.
\end{definition}

\begin{definition} {\rm(see \cite[Def. 2.1 and Def. 2.3]{chu-15})}
  {\rm \x{eq1}} admits   a   nonuniform  exponential  dichotomy if there exist an invariant projection  $P$, and constants $\alpha>0$, $M>0$, and $\varepsilon \in [0, \alpha)$ such that
 \begin{equation}\lb{b7}
\|\Phi(t)P\Phi^{-1}(s)\|\le M e^{-\alpha(t-s)}e^{\varepsilon s}, \quad {\rm for} \quad 0\le s \le t,
\end{equation}
and
\begin{equation}\lb{b8}
\|\Phi(t)Q\Phi^{-1}(s)\|\le M e^{\alpha(t-s)}e^{\varepsilon s}, \quad {\rm for} \quad 0\le t \le s,
\end{equation}
where  $Q =I_{n}-P$ is  the   complementary   projection. Furthermore, for any fixed $\gamma\in \R$, write a   shifted   system
\begin{equation}\lb{b9}
  \dot{x}=\left[A(t)x-\gamma I_{n}\right]x.
\end{equation}
Then the   nonuniform  dichotomy  spectrum  of {\rm \x{eq1}} is   the   set
\[\Sigma_{NED}(A)=\{\gamma\in \R: ~{\rm \x{b9}}~admits~no~nonuniform~ exponential~ dichotomy\},\]
and   the   resolvent   set $\rho_{NED}(A)=\R \setminus \Sigma_{NED}(A)$ is its complements.
\end{definition}
From the definition above, we know that the nonuniformity  means that $M$ is  no longer  a constant in the definition of mean-square exponential dichotomy  but a function $M e^{\varepsilon |s|}$ depending on the initial time $s$  (see \x{b7} and \x{b8} for details). We emphasize that in comparison to the notion of (uniform) exponential dichotomies(\cite{ss, ss1}), this notion is a weaker requirement. In particular,
 when $\varepsilon=0$, we obtain the notion of (uniform) exponential dichotomy.

The existence of integral separateness  can be used to prove the necessary and sufficient conditions for the stability of system \x{eq1} varies
under small perturbations (see e.g., \cite[Thm. 5.4.7]{ad-95} or \cite{by-69}). However,
as mentioned above, it is increasing interest to look for more general types about necessary and sufficient conditions for the stability of system \x{eq1} due to the fact that integral separateness needs more dynamical information, just as
Example \ref{example1} shows that \x{b10} is not integral separateness but weak integral separateness.

Before the proof Theorem \ref{main1}, it is important to mention that Barabanov and Denisenko \cite{bd-07} also establish  a necessary and sufficient condition for the stability of the Lyapunov exponents of \x{eq1} with exponential decaying perturbations by using exponentially  integral separation and higher and lower Izobov exponential indices. Unlike the results in \cite{bd-07}, Theorem \ref{main1} extends the condition of the coefficients matrix from uniformly bounded to nonuniformly bounded growth and establish the connection between exponential decaying perturbations and weak integral separation by using a different proof.

\vspace{0.5cm}

\noindent\mbox{\bf Proof of Theorem \ref{main1}.~}{ 
\emph{(Sufficiency)} Firstly, we illustrate that weak integral separateness is invariant under Lyapunov transformation. For every $j=1,\dots,n$, let $e_{j}$ be the unit column-vector in the $x_{j}$ direction, i.e., \[e_{j}=(\underbrace{0,\ldots,0,1}_{j},0,\ldots,0)^{T}.\]
Let $\Phi_{j}(t)=(\Phi_{1j}(t),\ldots,\Phi_{nj}(t))$ be the solution of \x{eq1} with initial value $\Phi_{j}(0) =e_{j}$. Thus $\Phi(t)=(\Phi_{1}(t),\ldots,\Phi_{n}(t))$ is a principal matrix solution of \x{eq1}. Let a weakly integrally separated system \x{eq1}
be reduced to the system $\dot{y}=B(t)y$ under a Lyapunov transformation $y=T(t)x$. we show the fundamental matrix $T(t)\Phi(t)=(T(t)\Phi_{1}(t),\ldots,T(t)\Phi_{n}(t))$ of this system is also weakly integrally separated.

Notice that $T$ is a Lyapunov transformation, thus there exists a constant $L>0$ such that
\[\|T(t)\|\leq L, \quad \|T^{-1}(t)\|\leq L \quad {\rm for } \quad t\ge 0.\]
Thus it follows easily from the inequality
\[\|T^{-1}(t)T(t)\Phi(t)\|\leq L \|T(t)\Phi(t)\|\]
that \[\|T(t)\Phi(t)\| \geq \frac{\|\Phi(t)\|}{L}.\]
This inequality, under the condition $\x{b5}$,  combined with $\|T(t)\Phi(t)\|\leq L\|\Phi(t)\|$ imply that
\[\frac{\|T(t)\Phi_{i+1}(t)\|}{\|T(t)\Phi_{i+1}(s)\|} \cdot \frac{\|T(t)\Phi_{i}(s)\|}{\|T(t)\Phi_{i}(t)\|} \geq \frac{1}{L^{4}}
\frac{\|\Phi_{i+1}(t)\|}{\|\Phi_{i+1}(s)\|} \cdot \frac{\|\Phi_{i}(s)\|}{\|\Phi_{i}(t)\|}\geq \frac{D}{L^{4}} e^{a(t-s)-b s}
\]
with $D>0$ and $a, b \geq 0$.

Secondly, we show that a weakly
integrally separated system is reducible to a diagonal one by using the Lyapunov transformation. Assuming that the inequality \x{b5} is satisfied with the principal matrix solution $\Phi(t)=(\Phi_{1}(t),\ldots,\Phi_{n}(t))$ of system \x{eq1}. From Corollary 3.3.2 and Remark 3.3.4 in \cite{ad-95}, we know that any linear system can be reducible to a diagonal form by methods of a Lyapunov transformation if and only if  the fundamental matrix $\Phi(t)$ satisfies the condition
\[\frac{G(\Phi)}{\|\Phi_{1}(t)\|^{2}\cdot \|\Phi_{2}(t)\|^{2}\cdots \|\Phi_{n}(t)\|^{2}}=\sin^{2} \alpha_{1}\cdots \sin^{2} \alpha_{n-1}\ge \varrho >0 \quad {\rm for} \quad t\in \R_{+}\]
with $G(\Phi)$ is the Gram determinant of the fundamental matrix solution $\Phi(t)$, and
\[\alpha_{k}=\sphericalangle \left(L_{k}, \Phi_{k+1}(t)\right), \quad k=1,\ldots, n-1,\]
is the angle between $\Phi_{k+1}(t)$ and $L_{k}$, where $L_{k}$ is the $k-$dimensional vector subspace spanned by the solutions $\Phi_{1}(t),\Phi_{2}(t),\ldots, \Phi_{k}(t)$.

Notice that the condition of nonuniformly bounded growth \x{eq2} is satisfied, then  $\Sigma_{NED}(A)$ is a bounded   closed   set  and $\Sigma_{NED}(A) \in [-\tilde{a},\tilde{a}]$ (See e.g., \cite[Lemma 2.10]{chu-15}). Thus it follows from \x{eq2} that
\[\frac{\|\Phi_{i}(s)\|}{\|\Phi_{i}(t)\|}\ge \frac{1}{K} e^{-\tilde{a}(t-s)}e^{-\tilde{b} s},  \quad {\rm for} \quad 0\le s \le t,~1\le i \le n.\]
Hence, the shifted system
\[\dot{x}=(A(t)+\lambda I_{n} )x\]
of \x{eq1} satisfies
\be\lb{b11}\frac{\|\Phi^{-\lambda}_{i}(s)\|}{\|\Phi^{-\lambda}_{i}(t)\|}\ge \frac{1}{K} e^{(\lambda-\tilde{a})(t-s)}e^{-\tilde{b} s}, \quad {\rm for} \quad 0\le s \le t,~1\le i \le n,\ee
with $\lambda$ is sufficiently large such that $\lambda-\tilde{a}>0$. In fact, this transformation, increases the characteristic exponents of \x{eq1} but does not affect in any way the angles and the stability of system \x{eq1}. In order to simplify the presentation, we omit the shift coefficient index $\lambda$, i.e., let $\Phi_{i}(t)=\Phi^{-\lambda}_{i}(t)$. Then we try to use induction to show that all the angles $\alpha_{k}$ are bounded away from zero.

 It is trivial for $k=1$. Assuming that $\alpha_{k}$ are bounded away from zero for all $k=2,\ldots,m$, i.e.,
 \be\lb{b12}\frac{G(\Phi_{1},\ldots,\Phi_{k})}{\|\Phi_{1}\|^{2}\cdot \|\Phi_{2}\|^{2}\cdots \|\Phi_{k}\|^{2}}=\sin^{2} \alpha_{1}\cdots \sin^{2} \alpha_{k-1}\ge \varrho >0.\ee
Now we prove \x{b12} for $k=m+1$. Assume the contrary, that there exists a sequence of solutions $\tilde{\Phi}_{m}(t_{i})\in L_{m}$ such that
\be\lb{b13}\alpha_{m}(t_{i})=\sphericalangle \left(\tilde{\Phi}_{m}(t_{i}), \Phi_{m+1}(t_{i})\right)\rightarrow 0,\quad {\rm as } \quad t_{i} \rightarrow +\oo.\ee
Without loss of generality, we assume that
\[\|\tilde{\Phi}_{m}(s_{i})\|=\|\Phi_{m+1}(s_{i})\|=1.\]
Then,
\be\lb{b14}\|\Phi_{m+1}(t_{i})-\tilde{\Phi}_{m}(t_{i})\|\rightarrow 0 \quad {\rm as } \quad t_{i} \rightarrow +\oo.\ee

On the other side, it follows from \x{eq2}, \x{b5},  \x{b11} and $\lambda-\tilde{a}>0$ that
\beaa \EM \|\Phi_{m+1}(t_{i})-\tilde{\Phi}_{m}(t_{i})\|\\ \GE \|\Phi_{m+1}(t_{i})\|-\|\tilde{\Phi}_{m}(t_{i})\|\\
\EQ \|\Phi_{m+1}(s_{i})\|\frac{\|\Phi_{m+1}(t_{i})\|}{\|\Phi_{m+1}(s_{i})\|}
-\|\tilde{\Phi}_{m}(s_{i})\|\frac{\|\Phi_{m+1}(t_{i})\|}{\|\Phi_{m+1}(s_{i})\|}
\left(\frac{\|\tilde{\Phi}_{m}(t_{i})\|}{\|\tilde{\Phi}_{m}(s_{i})\|}
\frac{\|\Phi_{m+1}(s_{i})\|}{\|\Phi_{m+1}(t_{i})\|}\right)\\
\GE \frac{1}{K} e^{(\lambda-\tilde{a})(t_{i}-s_{i})}e^{-\tilde{b} s_{i}} \left(1- \frac{1}{D} e^{-a(t_{i}-s_{i})}e^{b s_{i}} \right)>\frac{1}{2}.\eeaa
The last inequality above can be guaranteed by letting $t_{i}-s_{i}$ sufficiently large for any fixed $s_{i}$. The inequality obtained above contradicts the condition \x{b14}, which is due to the hypothesis \x{b13}. Thus the inequality \x{b13} holds for $k=m+1$.

Now we prove that \x{eq1} is stable with the perturbations of the coefficient matrix being exponential decaying. Assume that diagonal system \x{eq1} is weakly integrally separated, from the above proof, there exists a Lyapunov transformation, which can reduce \x{eq1} into a diagonal one such that diagonal coefficients are ordered as
\[\int^{t}_{s}(a_{i+1}(\tau)-a_{i}(\tau)) d \tau\geq a(t-s)-b s - D, \quad i=1,\ldots,n-1,\quad t\ge s \ge 0\] with $a\ge 0$, $b \geq 0$ and $D\in \R$.
 It follows from the condition of nonuniformly bounded growth \x{eq2} that the Cauchy matrix of \x{b6} defined by $\Phi(t,s)=\Phi(t)\Phi^{-1}(s)$
satisfies
\[\|\Phi(t,s)\|\le K e^{\tilde{a}(t-s)}e^{\tilde{b} s}, \quad {\rm for} \quad 0\le s \le t\]
with $\Phi(t)$ is a principal matrix solution of \x{b6}. This mean that the shifted system
\[\dot{y}=(diag [a_{1}(t),\ldots, a_{n}(t)]-\lambda I_{n})y\]
of \x{b6} satisfies a nonuniform exponential contraction
\be\lb{b16}\|\Phi^{\lambda}(t,s)\|\le K e^{-(\lambda-\tilde{a})(t-s)}e^{\tilde{b} s}, \quad {\rm for} \quad 0\le s \le t\ee with $\lambda-\tilde{a}>0$ and $\tilde{b}\geq 0$. As mentioned above, this transformation does not affect in any way the stability of system \x{b6}. In order to simplify the presentation, we omit the shift coefficient index $\lambda$, i.e., we set $\Phi(t,s)=\Phi^{\lambda}(t,s)$.
To prove that the Lyapunov exponents of \x{b6} are stable it suffices to show that
 the Lyapunov exponents of the shifted system \x{b16} satisfy the condition of the stability.

Now let $\Psi(t,s)$ be the Cauchy matrix of the perturbed system
\be\lb{b17}\dot{y}=(diag [a_{1}(t),\ldots, a_{n}(t)]-\lambda I_{n}+B(t))y.\ee
Let \[J=\{(t,s)\in \R^{+} \times \R^{+}:~ t\ge s \ge 0\},\]
and set
\[\mathcal{X}=\{\Psi: J\rightarrow \mathcal{B}(\R^{n}): \Psi~ is~continuous~and ~\|\Psi\|<\oo \},\]
which is a Banach space with the norm
\[\left\|\Psi\right\|_{\mathcal{X}}=\sup\{\|\Psi(t,s)\|e^{-bs}:(t,s)\in J\}.\]
Let \[(L\Psi)(t,s)=\Phi(t,s)+\int_{s}^{t}\Phi(t,\tau)B(\tau)
\Psi(\tau,s) d \tau\] for every $\Psi\in \mathcal{X}$.  Then, for each $\Psi_{1}, \Psi_{2}\in \mathcal{X}$, one has
\beaa \|(L\Psi_{1})(t,s)-(L\Psi_{2})(t,s)\| \LE
\int_{s}^{t}\|\Phi(t,\tau)\|\cdot\|B(\tau)\|\cdot \|\Psi_{1}(\tau,s)-\Psi_{2}(\tau,s)\|d \tau\\
\LE K\delta e^{bs}\|\Psi_{1}-\Psi_{2}\|_{\mathcal{X}}\int_{s}^{t}e^{-(\lambda-\tilde{a})(t-\tau)} d \tau. \eeaa
This means that \[\|(L\Psi_{1})-(L\Psi_{2})\|_{\mathcal{X}}\le \frac{K \delta}{\lambda-\tilde{a}}\|\Psi_{1}-\Psi_{2}\|_{\mathcal{X}}\] due to the fact $\lambda-\tilde{a}>0$. Moreover, $\delta<(\lambda-\tilde{a})/K$ implies the operator $L: \mathcal{X} \rightarrow \mathcal{X}$ is a contraction. Hence, there exists a unique $\Psi\in \mathcal{X}$ such that $L\Psi=\Psi$ which satisfies
\[\Psi(t,s)=\Phi(t,s)+\int_{s}^{t}\Phi(t,\tau)B(\tau)
\Psi(\tau,s) d \tau.\]
This gives the inequality of the form
\be\lb{b18}\|\Psi(t,s)\|\le \|\Phi(t,s)\|+\int_{s}^{t}\|\Phi(t,\tau)\|\cdot \|B(\tau)\|\cdot
\|\Psi(\tau,s)\| d \tau.\ee
It follows from \x{b18} that
\[\left\|\frac{\Psi(t,s)}{\Phi(t,s)}\right\| \le
1+  \int_{s}^{t}\|B(\tau)\|\cdot\left\|\frac{\Psi(\tau,s)}{\Phi(\tau,s)}\right\| d \tau.
\]
Set now $z(t):=\|\Psi(t,s)/\Phi(t,s)\|$ with any fixed initial point $s\in[0,t]$, the inequality above implies that
\[z(t)\le e^{\int_{s}^{t}\|B(\tau)\| d \tau}\le e^{\delta (t-s)}.\]
By using the inequality above, and the equality
\[\limsup_{t\rightarrow \oo}\frac{1}{t}\ln\|f(t)\|=-\liminf_{t\rightarrow \oo}\frac{1}{t}\ln\|1/f(t)\|,\]
one can easily verify that
\beaa|\lambda_{i}-\hat{\lambda}_{i}|\EQ
|\lambda(\Psi(t,s)e_{i})-\lambda(\Phi(t,s)e_{i})|\\
\EQ \left|\limsup_{t\rightarrow \oo}\frac{\ln\|\Psi(t,s)e_{i}\|}{t}-\limsup_{t\rightarrow \oo}\frac{\ln\|\Phi(t,s)e_{i}\|}{t}\right|\\
\LE \left|\limsup_{t\rightarrow \oo}\frac{\ln\|\Psi(t,s)e_{i}\|}{t}-\liminf_{t\rightarrow \oo}\frac{\ln\|\Phi(t,s)e_{i}\|}{t}\right|\\
\EQ \left|\limsup_{t\rightarrow \oo}\frac{\ln\|\Psi(t,s)e_{i}\|}{t}+\limsup_{t\rightarrow \oo}\frac{\ln\|1/(\Phi(t,s)e_{i})\|}{t}\right|\\
\LE \delta\eeaa
for $i=1,\cdots,n$. Now the stability \x{b1} follows with $\delta<\varepsilon$.

\vspace{0.5cm}
\emph{(Necessity)} Assume that the system \x{b6}, which has distinct characteristic exponents $\lambda_{1} > \cdots > \lambda_{n}$, is stable with the perturbations of the coefficient matrix being exponential decaying, i.e., $\|B(t)\|\le \delta e^{-\beta t}$. By virtue of  the method of variation of constants, any  Cauchy matrix of
\be\lb{fj2}\dot{y}=(diag [a_{1}(t),\ldots, a_{n}(t)]e^{\beta t}+B(t)e^{\beta t})y.\ee
satisfies the integral equation
\[y(t)=\Phi(t,s)y(t_{0})+\int_{s}^{t}\Phi(t,\tau)B(\tau)e^{\beta \tau}
y(\tau) d \tau,\]
where  $\Phi(t,s)$ is the Cauchy matrix of
  \be\lb{fj1}\dot{y}=diag [a_{1}(t),\ldots, a_{n}(t)]e^{\beta t}y.\ee
 Note that $\|B(t)e^{\beta t}\| \le \delta$ for all $t\ge 0$,  following the same methods as in the proof of the sufficiency, we can prove that \x{fj1}
is stable with the perturbation $B(t)e^{\beta t}$.
Thus it follows from Bylov and Izobov's result \cite{by-69} that the diagonal elements of \x{fj1} are integrally separated, that is, there exist some constants $\hat{a}>0$ and $D\in \R$ such that
\be\lb{fj4}
\int^{t}_{s}(a_{i+1}(\tau)-a_{i}(\tau))e^{\beta \tau} d \tau\geq \hat{a}(t-s)- D, \quad t\ge s \ge 0.
\ee
Thus, to prove \x{b6} is weakly integrally separated, it suffices to prove that
\be\lb{fj5}
\int^{t}_{s}(a_{i+1}(\tau)-a_{i}(\tau)) d \tau\geq a(t-s)-bs- D, \quad t\ge s \ge 0
\ee
for $i=1, \ldots, n-1$ and some constants $a, b \ge 0$ and $D\in \R$.

In fact, if $a_{i+1}(\tau)-a_{i}(\tau)\ge 0$ on any interval $[s,t]$, the inequality \x{fj5} holds with $a, b, D =0$. Conversely, if $a_{i+1}(\tau)-a_{i}(\tau)\le 0$ on any interval $[s,t]$, the condition \x{fj4} does not hold due to the fact that the right side of \x{fj4} is always positive with any interval $[s,t]$ large enough.

Now we prove the nontrivial case: $a_{i+1}(\tau)-a_{i}(\tau)\le 0$ on the disjoint union of finite closed intervals, i.e., $\bigcup_{k}[s_{k},t_{k}]$. It follows from \x{fj4} that
\beaa D\GE \int^{t_{k}}_{s_{k}}-(a_{i+1}(\tau)-a_{i}(\tau))e^{\beta \tau} d \tau+\hat{a}(t-s)\\
\GE \int^{t_{k}}_{s_{k}}-(a_{i+1}(\tau)-a_{i}(\tau)) d \tau+\hat{a}(t-s).\eeaa
This gives \x{fj5} with $a=\hat{a}$, and $b=0$. Thus the proof of the necessity is complete. \hspace{\stretch{1}}$\Box$
\\

\vspace{0.2cm}
Let $\tilde{z}_{1}(t):=\|\Psi(t,s)\|$ with any fixed initial point $s\in[0,t]$, it follows from  \x{b16} and \x{b18} that for $ \|B(t)\|\le \delta e ^{-b s}$,
\[\tilde{z}_{1}(t)\le K e^{-(\lambda-\tilde{a})(t-s)+b s} + \delta K \int_{s}^{t}e^{-(\lambda-\tilde{a})(t-\tau)}\tilde{z}_{1}(\tau)d \tau \quad t\ge s \ge 0.\]
Consider the continuous function $\tilde{z}_{2}(t)$ satisfies the integral function
\be\lb{b19} \tilde{z}_{2}(t)= K e^{-(\lambda-\tilde{a})(t-s)+b s} + \delta K \int_{s}^{t}e^{-(\lambda-\tilde{a})(t-\tau)}\tilde{z}_{2}(\tau)d \tau \quad t\ge s \ge 0.\ee It is easy to prove that the integral function \x{b19} is equivalent to the differential equation $\tilde{z}'_{2}(t)=(\delta k -\lambda +\tilde{a})\tilde{z}_{2}(t)$ with the initial condition $\tilde{z}_{2}(s)=K e^{bs}$. Hence, there exist a unique solution \[\tilde{z}_{2}(t)=K e^{(\delta k -\lambda +\tilde{a})(t-s)+b s}\] of the integral equation \x{b19}. Thus, obviously,
\[\|\Psi(t,s)\|\le \tilde{z}_{2}(t)=K e^{-(\lambda -\tilde{a}-\delta k )(t-s)+b s}.\]
This shows that the Cauchy matrix $\Psi(t,s)$ of the perturbed system \x{b17} also admits a nonuniform exponential contraction with $\delta>0$ small enough. Furthermore, Cauchy matrix $\Psi(t,s)$ of the perturbed system \x{eq3} admits a nonuniform exponential contraction follows from the proof the Theorem \ref{main1} with $\tilde{a}<0$. More precisely, the trivial solution of  \x{eq3} is also asymptotically and exponentially stable  with sufficiently small $\delta>0$ if the trivial solution of  \x{eq1} is asymptotically and exponentially stable.

\begin{corollary}\lb{main2} Assume that the system {\rm\x{eq1}} with nonuniformly bounded growth  has distinct Lyapunov exponents $\lambda_{1} > \cdots > \lambda_{n}$.
Then they are stable with the perturbations of the coefficient matrix being exponentially decaying if and only if there exists a fundamental matrix solution with weakly integrally separated columns.
\end{corollary}
\noindent\prf{For the proof of Corollary \ref{main2},  it suffices to show that a fundamental matrix solution of \x{eq1}  with weakly integrally separated columns and
the weakly integrally separated diagonal coefficients $a_{i}$ of \x{b6} are
equivalent. In fact, assume that $\Phi(t)=(\Phi_{1}(t),\ldots,\Phi_{n}(t))$ is a fundamental matrix solution of system \x{eq1}, and let
\[T^{-1}(t)=\left\{\frac{\Phi_{1}(t)}{\|\Phi_{1}(t)\|},
\frac{\Phi_{2}(t)}{\|\Phi_{2}(t)\|},\cdots,\frac{\Phi_{n}(t)}{\|\Phi_{n}(t)\|}\right\}\]
 be a Lyapunov transformation, which satisfies
 \[T(t)\Phi(t)=diag[\|\Phi_{1}(t)\|,\ldots,\|\Phi_{n}(t)\|].\]
This implies \[Y(t)=diag[\|\Phi_{1}(t)\|,\ldots,\|\Phi_{n}(t)\|]\] is the fundamental matrix of \x{b6} and
\beaa diag [a_{1}(t),\ldots, a_{n}(t)]\EQ \frac{\dot{Y}(t)}{Y(t)}=diag\left[\frac{d}{dt}\ln\|\Phi_{1}(t)\|,\ldots,
\frac{d}{dt}\ln\|\Phi_{n}(t)\|\right].\eeaa
Then we have
\beaa\int^{t}_{s}(a_{i+1}(\tau)-a_{i}(\tau)) d \tau \EQ \ln \left(\frac{\|\Phi_{i+1}(t)\|}{\|\Phi_{i+1}(s)\|} \cdot \frac{\|\Phi_{i}(s)\|}{\|\Phi_{i}(t)\|}\right),\quad  i=1, \ldots, n-1,
\eeaa
for $0\le s\le t$, and this means that
\[\frac{\|\Phi_{i+1}(t)\|}{\|\Phi_{i+1}(s)\|} \cdot \frac{\|\Phi_{i}(s)\|}{\|\Phi_{i}(t)\|}\geq e^{D} e^{a(t-s)-b s}
\Leftrightarrow
\int^{t}_{s}(a_{i+1}(\tau)-a_{i}(\tau)) d \tau\geq a(t-s)-b s + D\]
with $a, b \geq 0$ and $D\in \R$.
\hspace{\stretch{1}}$\Box$
}

\section{\bf{Weak integral separateness and nonuniform exponential dichotomy spectrum}}
\setcounter{equation}{0} \noindent

In this section we try to establish the connection between weak integral separateness and nonuniform exponential dichotomy spectrum.

The following theorem can be used to illustrate the existence of weak integral separateness under full spectrum $\Sigma_{NED}(A)$. This result of uniform type ($b=0$ in \x{b5}) is given in  \cite[p. 231]{bs-00}.  Here we give a different proof based on the shift system.

\begin{theorem}\lb{main3}
  Assume that {\rm\x{eq1}} has a  full  nonuniform spectrum, i.e., \[\Sigma_{NED}(A)=\bigcup_{i=1}^{n}[a_{i},b_{i}].\]
  Then there exists a fundamental matrix solution with weakly integrally separated columns.
\end{theorem}
\noindent\prf{ Let $\lambda_{1}=(a_{2}+b_{1})/2$. Obviously, $\lambda_{1}\in \rho_{NED}$, and this means that a fundamental matrix $\Phi_{\lambda_{1}}(t)$ of the shift system
\[\dot{x}=\left[A(t)x-\lambda_{1} I_{n}\right]x \]
admits a nonuniform exponential dichotomy
\beaa
\|\Phi_{\lambda_{1}}(t)P\Phi_{\lambda_{1}}^{-1}(s)\|\LE M e^{-\alpha(t-s)}e^{\varepsilon s}, \quad {\rm for} \quad 0\le s \le t, \\
\|\Phi_{\lambda_{1}}(t)(I_{n}-P)\Phi_{\lambda_{1}}^{-1}(s)\|\LE M e^{\alpha(t-s)}e^{\varepsilon s}, ~\quad {\rm ~~for} \quad 0\le t \le s,
\eeaa
with $\alpha>0$, $\varepsilon\in [0,\alpha)$, and
\[P=\begin{pmatrix}
      I_{(n-1) \times (n-1)} & 0_{(n-1) \times 1} \\
      0_{1 \times (n-1)} & 0 \\
    \end{pmatrix}
\]
Thus, following the method in \cite[Lemma 5.1.1]{ad-95}, we can prove
\bea\lb{b20}\frac{\|\Phi_{1}(t)\|}{\|\Phi_{1}(s)\|}\cdot e^{\lambda_{1}(t-s)}\EQ
\|\Phi_{\lambda_{1}}(t)P\Phi_{\lambda_{1}}^{-1}(s)\| \le M e^{-\alpha(t-s)}e^{\varepsilon s}, \quad {\rm for} \quad 0\le s \le t,\eea
and
\bea\lb{b21}\frac{\|\Phi_{j}(t)\|}{\|\Phi_{j}(s)\|}\cdot e^{\lambda_{1}(t-s)}\LE
\max_{\|b\|=1}\|\Phi_{\lambda_{1}}(t)(I_{n}-P)\Phi_{\lambda_{1}}^{-1}(s)b\|\le M e^{\alpha(t-s)}e^{\varepsilon s}\eea
for $0\le t\le s$, and $j=2,\ldots, n$.
Then
\[\frac{\|\Phi_{2}(t)\|}{\|\Phi_{2}(s)\|}\cdot \frac{\|\Phi_{1}(s)\|}{\|\Phi_{1}(t)\|}\ge \frac{1}{M^{2}}e^{2\alpha(t-s)}e^{2\varepsilon s}.\]
 Repeating the procedure \x{b20}-\x{b21} for $\lambda_{i}=(a_{i+1}+b_{i})/2$ for all $i=2,\ldots , n-1$, and this completes the proof. \hspace{\stretch{1}}$\Box$}

\vspace{0.5cm}
  On the contrary, the result is not true. In the  following, we provides a simple example illustrating that  weak integral separateness does not guarantee the existence of
  full  nonuniform dichotomy spectrum.

 \begin{example}\lb{example31}
 Consider a diagonal system
 \be\lb{b22}\dot{x}_{1}=(2-2t\sin t)x_{1},\quad \dot{x}_{2}=(4-3t\sin t)x_{2}.\ee
 From {\rm \cite[Example 2.1]{chu-15}}, we know that $\Sigma_{NED}(A)=[\lambda-a,\lambda+a]$ for \[\dot{x}=(\lambda-a t\sin t)x.\]
 Hence, the nonuniform dichotomy spectrum $\Sigma_{NED}=[0,4]\cup[1,7]$ of {\rm \x{b22}} overlap and interact. Meanwhile,
 \[\frac{|x_{2}(t)|}{|x_{2}(s)|}\cdot \frac{|x_{1}(s)|}{|x_{1}(t)|}\ge e^{t-s}e^{2s},\quad 0\le s \le t,\]
 so that $x_{1}$ and $x_{2}$ are weakly integrally separated. \hspace{\stretch{1}}$\Box$
 \end{example}

Based on this research, we have to find more information about full  nonuniform dichotomy spectrum.  In order to achieve this goal, we first consider the existence conditions of nonuniform exponential dichotomy.
In fact, in the case of constant coefficients, only  uniform exponential dichotomy can exist if and only if the eigenvalues of the coefficient
matrix have nonzero real parts. In view of this idea,  if one wants to prove that a  linear differential system has a nonuniform exponential dichotomy, it is necessary to find a  subspace of solutions with  a nonuniformly bounded growth and a complementary subspace  of solutions with  a nonuniformly bounded decay.  the following theorem establish a sufficient but not necessary condition for the existence of a nonuniform exponential dichotomy with $A(t)$  bounded away from zero.

\begin{theorem}\lb{main4}Assume that system {\rm\x{eq1}} with nonuniformly bounded growth  is weakly integrally separated, and $\inf_{t\ge 0}|\det A(t)|>0$ holds. Then {\rm\x{eq1}} has a nonuniform exponential dichotomy {\rm\x{b7}-\x{b8}}.
\end{theorem}
\noindent\prf{
 It is easy to see from the proof of Theorem \ref{main1} that a weakly
integrally separated system is reducible to a diagonal one by using the Lyapunov transformation. Moreover, nonuniform exponential dichotomy is invariant under Lyapunov transformation. In fact, assume that \x{b6} admits a nonuniform exponential dichotomy for a fundamental matrix $\tilde{\Phi}(t)$ with a Lyapunov transformation $y\leftarrow T(t)x$ such that $T(t)\tilde{\Phi}(t)=\Phi(t)$, where $\Phi(t)$ is a fundamental matrix solution of \x{eq1}.  Then
by the condition of Lyapunov transformation, i.e., there exists a constant $L>0$ such that
\[\|T(t)\|\leq L, \quad \|T^{-1}(t)\|\leq L \quad {\rm for } \quad t\ge 0,\]
we obtain
\beaa \|\Phi(t)P\Phi^{-1}(s)\| \EQ \|T(t)\tilde{\Phi}(t)P\tilde{\Phi}^{-1}(s) T^{-1}(s)\|\\
\LE \|T(t)\|\cdot \|\tilde{\Phi}(t)P\tilde{\Phi}^{-1}(s)\|\cdot \|T^{-1}(s)\| \\
\LE L^2M e^{-\alpha(t-s)}e^{\varepsilon s},\quad 0<s\le t.
\eeaa
A   similar   argument  shows   that
\beaa \|\Phi(t)(I_{n}-P)\Phi^{-1}(s)\|
\LE L^2M e^{\alpha(t-s)}e^{\varepsilon s},\quad 0<t\le s,
\eeaa
and hence the fundamental matrix solution $\Phi(t)$ admits a nonuniform exponential dichotomy. Thus, it suffices to prove that \x{b6} with nonuniformly bounded growth admits a nonuniform exponential dichotomy.

Naturally, since \x{eq1} is weakly integrally separated, it follows from Theorem \ref{main1} that system \x{b6} can be rewritten as
\beaa
\dot{x}_{i}\EQ a_{i}(t)x_{i} \quad (i=1,\ldots, k-1)\\
\dot{x}_{k}\EQ a_{k}(t)x_{k}\\
\dot{x}_{i}\EQ a_{i}(t)x_{i} \quad (i=k+1,\ldots, n),
\eeaa
with
\be\lb{b23} exp\left(\int_{s}^{t}(a_{i+1}(\tau)-a_{i}(\tau))d \tau \right)\geq D e^{a(t-s)-b s}, \quad i=1,\ldots, n-1, \quad t\ge s\ge 0.\ee

Note that $\inf_{t\ge 0}|\det A(t)|>0$,  We assume, without loss of generality, that
$a_{k}(t)<-\epsilon<0$. Thus the scalar equation
\[\dot{x}_{k}= a_{k}(t)x_{k}\] has an exponential dichotomy, i.e., for $\epsilon>0$ such that for $0\le s\le t$,
\[exp\left(\int_{s}^{t}a_{k}(\tau)d \tau \right)\leq e^{-\epsilon(t-s)}.\]
Then, by \x{b23},
\beaa exp\left(\int_{s}^{t}a_{k-1}(\tau)d \tau \right) \EQ exp\left(\int_{s}^{t}(a_{k-1}(\tau)-a_{k}(\tau))d \tau \right)\cdot exp\left(\int_{s}^{t}a_{k}(\tau)d \tau \right)\\
\LE \frac{1}{D}e^{-a(t-s)+b s}e^{-\epsilon(s-t)}=e^{-(a+\epsilon)(t-s)+b s}.
\eeaa
Repeating this argument we see that each of the first k equations in \x{b6} has a nonuniform exponential dichotomy. If $k=n$ we are finished with $P=I_{n}$. So suppose that $k<n$, thus it follows from $\inf_{t\ge 0}|\det A(t)|>0$ that $a_{k+1}(t)>\epsilon>0$. Thus the scalar equation
\[\dot{x}_{k+1}= a_{k+1}(t)x_{k+1}\] has an exponential dichotomy, i.e., for $\epsilon>0$ such that for $0\le t\le s$,
\[exp\left(\int_{s}^{t}a_{k+1}(\tau)d \tau \right)\leq e^{-\epsilon(s-t)}.\]
Then, by \x{b23},
\beaa exp\left(\int_{s}^{t}a_{k+2}(\tau)d \tau \right) \EQ exp\left(\int_{s}^{t}(a_{k+2}(\tau)-a_{k+1}(\tau))d \tau \right)\cdot exp\left(\int_{s}^{t}a_{k+1}(\tau)d \tau \right)\\
\LE \frac{1}{D}e^{-a(t-s)+b s}e^{-\epsilon(s-t)}=e^{-(a+\epsilon)(t-s)+b s}.
\eeaa
Repeating this argument we can deduce each of the last $(n-k)$ equations in \x{b6} has a nonuniform exponential dichotomy. This means that \x{b6}, and then \x{eq1} has a nonuniform exponential dichotomy, and thus we are finished.  \hspace{\stretch{1}}$\Box$
}

\vspace{0.5cm}
Based on the observation of nonuniformly bounded growth, the aim of the following theorem is to give a sufficient condition for the existence of nonuniform dichotomy spectrum.

\begin{theorem}\lb{main5}
  Assume that system {\rm\x{b6}} with nonuniformly bounded growth is weakly integrally separated. Considering a matrix solution $\Phi(t)=(\Phi_{1}(t),\ldots,\Phi_{n}(t))$ of {\rm\x{b6}} such that distinct Lyapunov exponents of the columns of $\Phi$ are ordered as $\lambda_{1} > \cdots > \lambda_{s}$ for some $s\le n$ with some multiplicity $k_{i}$ at $\lambda_{i}$, i.e., $\sum_{m=1}^{s}k_{m}=n$. If there exists a $\lambda \in (\lambda_{j}, \lambda_{j+1})$ such that for $j=\sum_{m=1}^{i}k_{m}$,
  \be\lb{b24}\dot{y}=\left(
       \begin{array}{ccccccc}
         a_{1}(t) &   &  &   &   &   &   \\
           & \ddots &   &  &   &   &  \\
           &   & a_{j}(t)&  &  &   &  \\
          &  &  & \lambda &  &  &  \\
           &  &   &  & a_{j+1}(t) &   &   \\
           &  &  &   &   & \ddots &   \\
           &   &   &   &   &   &  a_{n}(t) \\
       \end{array}
     \right)y
  \ee
  is weakly integrally separated. Then the nonuniform dichotomy spectral  intervals can be splitted into two disconnected  subdomains.
\end{theorem}
  \noindent\prf{
  Firstly, it is trivial to know that Lyapunov exponents is contained within nonuniform dichotomy spectral  intervals. In fact, let $\lambda_{j}=\lambda(\Phi(t)P_{j}\Phi^{-1}(s))\in \Sigma_{NED}(A)$, $j=1,\ldots,s$ be a Lyapunov exponent with any fixed initial point $s\in[0,t]$, where $P_{j}$ is a projection of the form
  $P_{j}=\left(
           \begin{array}{cc}
             I_{j} & 0 \\
             0 & 0 \\
           \end{array}
         \right)$ with $j=\sum_{m=1}^{i}k_{m}$.
Then the Lyapunov exponent of the shifted system
         \[\dot{y}=(diag [a_{1}(t),\ldots, a_{n}(t)]-\lambda_{j} I_{n})y\]
of \x{b6} is written as
  \[\lambda(\Phi_{\lambda_{j}}(t)P_{j}\Phi_{\lambda_{j}}^{-1}(s))=\limsup_{t\rightarrow \oo}\frac{\|\Phi_{\lambda_{j}}(t)P_{j}\Phi_{\lambda_{j}}^{-1}(s)\|}{t}=0,\] which contradicts to the condition
  \beaa
\|\Phi_{\lambda_{j}}(t)P_{j}\Phi_{\lambda_{j}}^{-1}(s)\|\LE M e^{-\alpha(t-s)}e^{\varepsilon s}, \quad {\rm for} \quad 0\le s \le t, \\
\|\Phi_{\lambda_{j}}(t)(I_{n}-P_{j})\Phi_{\lambda_{j}}^{-1}(s)\|\LE M e^{\alpha(t-s)}e^{\varepsilon s}, ~\quad {\rm ~~for} \quad 0\le t \le s,
\eeaa
with $\alpha>0$, $\varepsilon\ge 0$. That is to say, $\lambda_{j} \in \rho_{NED}(A)$.

Now we show that the nonuniform dichotomy spectral  intervals can be splitted into two disconnected  subdomains. Thus, it follows from \x{b24} are weakly integrally separated and Lyapunov exponents of the columns of $\Phi$ are ordered as $\lambda_{1} > \cdots > \lambda_{s}$ that
\be\lb{b25} e^{\lambda(t-s)}e^{-\int_{s}^{t}a_{j}(\tau)d \tau}\ge D e^{a(t-s)-b s},  {\rm for} \quad 0\le s \le t,\ee
and
\be\lb{b26} e^{-\lambda(t-s)}e^{\int_{s}^{t}a_{j+1}(\tau)d \tau}\ge D e^{a(t-s)-b s}, {\rm for} \quad 0\le s \le t.\ee

From \x{b25}-\x{b26} we know that
\[e^{\int_{s}^{t}a_{j}(\tau)d \tau}\le \frac{1}{D}e^{-(\lambda+a)(t-s)+b s}, {\rm for} \quad 0\le s \le t,\]
and
\[e^{\int_{s}^{t}a_{j+1}(\tau)d \tau}\le \frac{1}{D}e^{-(\lambda-a)(t-s)+b s}, {\rm for} \quad 0\le t \le s.\]
Then following the same induction method in the proof of Theorem \ref{main4}, we can prove that each of the first $j$ equations and each of the last $(n-j)$ equations in the shifted system of \x{b6} has a nonuniform exponential dichotomy. This means that the nonuniform dichotomy spectral  intervals can be splitted into two disconnected  subdomains.  } \hspace{\stretch{1}}$\Box$


\vspace{0.5cm}

\noindent\mbox{\bf Proof of Theorem \ref{main11}.~}{By the proof of Theorem \ref{main1}, \x{eq1} is reducible to a diagonal system. Corollary 2.11 in  \cite{chu-15}  and Theorem \ref{main5} now imply the existence of the full nonuniform dichotomy spectrum.  } \hspace{\stretch{1}}$\Box$


}

\end{document}